\documentclass[reqno,11pt]{amsart}
\usepackage[cp1251]{inputenc}
\usepackage{graphicx}
\usepackage[russian]{babel}
\usepackage{amsmath,amssymb}
 \usepackage{eepic}
\usepackage{amsmath,amssymb}

\begin{document}
 \Large

\bigskip

\begin{center}

{\bf Hausdorff--Lebesgue dimension of attractors}

\bigskip

G.A.Leonov

\end{center}
\bigskip

\section{Introduction. Hausdorff measure and dimension and
Hausdorff--Lebesgue measure
and dimension}

In the present paper the classical ideas of Hausdorff and Lebesgue are combimed
and the Hausdorff--Lebesgue measure is introduced.
This makes it possible to  obtain new results in chaotic dynamics.

Consider a compact $K\subset R^n$ and the numbers $d\ge 0$, $\varepsilon>0$.

Define the Hausdorff measure and Hausdorff dimension of a compact $K$ [1,2].

Consider all coverings of $K$ by the balls $B_i$ of radii $r_i\le\varepsilon$.

Suppose,
$$
\mu_H(K,d,\varepsilon)=\inf\sum\limits_ir_i^d,
$$
where infimum is taken over all $\varepsilon$-coverings of compact $K$.

Obviously, $\mu_H(K,d,\varepsilon)$ increases with increasing
$\varepsilon$. Therefore there exists a limit
$$
\mu_H(K,d)=\lim\limits_{\varepsilon\to 0}\mu_H(K,d,\varepsilon).
$$

\textbf{Definition 1.} The value $\mu_H(K,d)$ is called a Hausdorff measure
of compact $K$.

We introduce
$$
\dim_HK=\inf\{d\,|\,\mu_H(K,d)=0\}.
$$

\textbf{Definition 2.} The value $\dim_HK$ is called a Hausdorff measure of $K$.

Note that a set of balls $\{B_i\}$ can be chosen as a set of cubes
with sides
$2r_i\le 2\varepsilon$. In this case the dimensions $\dim_HK$ coincide.

If the covering $\{B_i\}$ involves the balls of equal radii $r_i=\delta\le\varepsilon$,
we say about
fractal measure $\mu_F(K,d)$ and fractal
dimension $\dim_FK$.

The Hausdorff measure and fractal measure are outer measures.
But  in my view this measure is also outer.
Therefore in this paper we  combine the ideas of Hausdorff and Lebesgue.

Consider all coverings of $K$ by disjoint cubes $C_i$
with sides $2\delta_i\le 2\varepsilon$.

Also, as in the theory of Lebesgue measure,
in the case of intersection of boundaries $\partial C_i\cap\partial C_j$
such a set of intersections is included only in $C_i$ or in $C_j$.

Suppose that
$$
\mu_{HL}(K,d,\varepsilon)=\inf\sum\limits_I\delta_i^d,
$$
where the infimum is taken over all $2\varepsilon$-coverings of compact $K$.
It is obvious that $\mu_{HL}(K,d,\varepsilon)$ increases
with decreasing $\varepsilon$. Consequently there exists the limit
$$
\mu_{HL}(K,d)=\lim\limits_{\varepsilon\to 0}\mu_{HL}(K,d,\varepsilon).
$$

\textbf{Definition 3.} The value $\mu_{HL}(K,d)$
is called a Hausdorff--Lebesgue measure of compact $K$.

We introduce
$$
\dim_{HL}K=\inf\{d\,|\,\mu_{HL}(K,d)=0\}.
$$

\textbf{Definition 4.} The value $\dim_{HL}K$ is called a
Hausdorff--Lebesgue dimension of compact $K$.

Consider now all coverings of $K$ by disjoint cubes
$C_i$ with sides $2\varepsilon$

\textbf{Definition 5.} The value
$$
\mu_{FHL}(K,d)=\limsup\limits_{\varepsilon\to 0}\sum\limits_i\varepsilon^d
$$
is called a Hausdorff--Lebesgue fractal measure of compact $K$.

\textbf{Definition 6.} The value
$$
\dim_{FHL}K=\inf\{d\,|\,\mu_{FHL}(K,d)=0\}
$$
is called a Hausdorff--Lebesgue fractal dimension.

It is obvious that
$$
\mu_{H}(K,d,\sqrt n\varepsilon)\le\mu_{HL}(K,d,\varepsilon)(n)^{d/2},
$$
$$
\mu_{H}(K,d)\le \mu_{HL}(K,d)(n)^{d/2}\le\mu_{FHL}(K,d)(n)^{d\!/2},
$$
$$
\dim_HK\le\dim_{HL}K\le\dim_{FHL}K,
$$
and for $k$-dimensional manifold $K$
$$
\dim_HK=\dim_{HL}K=\dim_{FHL}K=k.
$$

The following relations are also obvious.

For compacts $K_i$ such that $K\subset\bigcup\limits_iK_i$
the inequality
$$
\mu_{HL}(K\bigcap(\bigcup\limits_i K_i),d,\varepsilon)
\le\sum\limits_i\mu_{HL}(K\bigcap K_i,d,\varepsilon)\eqno(1)
$$
is satisfied.
For disjoint compacts $K_i$ such
that $K\supset\bigcup\limits_i K_i$, the inequality
$$
\mu_{HL}(K\bigcap(\bigcup\limits_i K_i)d,\varepsilon)
\ge\sum\limits_i\mu_{HL}(K\bigcap K_i,d,\varepsilon)\eqno(2)
$$
is satisfied.
Similar relations are satisfied for $\mu_{FHL}(K,d)$.

\section{Upper Estimates of Hausdorff--Lebesgue dimension}

Recall [3] that a linear operator $A$ can be represented  in the form of
a product $A=SQ$,
where $S$ is symmetric nonnegative and $Q$ are orthogonal operators.
Recall also [3] that $S$ always has in $R^n$ orthonormal system of eigenvectors
$e_i $ $(j=1,\ldots,n)$ with real characteristic numbers that
coincide with singular values $\alpha_j$ $(\alpha_1\ge \cdots \ge\alpha_n\ge 0)$
of operator $A$.

\textbf{Definition 7.} A cube $C$ is called oriented
if the sides of cube $QC$ are parallel vectors $e_1,\ldots,e_n$.

Consider now a continuously differentiable mapping $F(x)\,:\,R^n\to R^n$
$$
F(x+h)-F(x)=(T_xF)h+o(h).\eqno(3)
$$
Suppose that $\alpha_1(T_xF)\ge\cdots\ge\alpha_n(T_xF)$ are
singular values of matrices $T_xF$ at the point $x$,
$$
\omega_d(T_xF)=\alpha_1(T_xF)\cdots\alpha_k(T_xF)\alpha_{k+1}(T_xF)^s,\quad d=k+s
$$

\textbf{Theorem 1.} Suppose that $FK=K$ and
$$
\sup\limits_K\omega_d(T_xF)<1.\eqno(4)
$$
Then
$$
\dim_{HL}K\le d.\eqno(5)
$$

Proof. From condition (4) it follows the existence of a number $\nu<1$ such that
$$
\sup\limits_K\omega_d(T_xF)\le \nu.\eqno(6)
$$
It is well known [2] that for a natural number $p$ it is valid the inequality
$$
\sup\limits_K\omega_d(T_xF^p)\le\nu^p.\eqno(7)
$$
We introduce the denotation
$$
\beta=\sup\limits_K\alpha_1(T_xF)
$$
$$
\gamma=\sup\limits_K\alpha_{k+1}(T_xF).
$$
It is obvious that $\gamma^d\le\nu$,
$$
\sup\limits_K\alpha_1(T_xF^p)\le\beta^p,\quad\sup\limits_K\alpha_{k+1}(T_xF^p)\le\nu^{p/d}.
$$
Choose $p$ in such a way that
$$
\sqrt n\nu^{p/d}\le 1,\quad 2^kn^{d/2}\nu^p<2^{-2}
$$
and $\varepsilon$ such that in the $\beta^p\varepsilon$-neighborhoods of
all points of compact $K$ there exists linearization procedure (3).

Consider a covering $K$ by the cubes $C_i$ with  sides $2\delta_i$
and centers at the points $x_i$. Consider also the oriented with respect to
 $T_{x_i}F^p$ cubes $\tilde C_i\supset C_i$ with centers $x_i$
and sides $2\sqrt n\delta_i$.

Obviously, $T_{x_i}F^p\tilde C_i$ is a
parallelipiped with sides $2\sqrt n\delta_i\alpha_j(T_{x_i}F^p)$, $j=1,\ldots,n$, and
$$
T_{x_i}F^pC_i\subset T_{x_i}F^p\tilde C_i.
$$
We cover this parallelipiped by cubes with sides
$2\sqrt n\delta_i\alpha_{k+1}(T_{x_i}F^p)$. The number of such cubes
is  less than or equal to
$$
\left(\frac{\alpha_1(T_{x_i}F^p)}{\alpha_{k+1}(T_{x_i}F^p)}+1\right)\cdots
\left(\frac{\alpha_k(T_{x_i}F^p)}{\alpha_{k+1}(T_{x_i}F^p)}+1\right).
$$
Consequently
$$
\mu_{HL}(T_{x_i}F^p\tilde C_i,d,\varepsilon)\le2^kn^{d/2}
\omega_d(T_{x_i}F^p)\delta_i^d\le\frac{\delta_i}{4}.
$$
Then by (1) we have
$$
\mu_{HL}(K,d,\varepsilon)=\mu_{HL}(F^pK,d,\varepsilon)<\frac{1}{2}\mu_{HL}(K,d,\varepsilon).
$$
However in this case $\mu_{HL}(K,d,\varepsilon)=0$ and $\mu_{HL}(K,d)=0$.
This implies the assertion of theorem.

\textbf{Theorem 2.} Suppose that for the compacts $\tilde K\supset K$
it is valid the following conditions $F^m(K)\subset\tilde K$, $\forall m\ge 1$,
$$
\sup\limits_{\tilde K}\omega_d(T_xF)<1,
$$
$$
\mu_{HL}(K,d)<\infty.
$$
Then
$$
\lim\limits_{m\to\infty}\mu_{HL}(F^m(K),d)=0.
$$
This theorem is an analog of the theorems on Hausdorff measure, proved in [4].
The proof of Theorem 2 is similar to the scheme, used in [4] with applying
the estimates, obtained in proving Theorem 1.

The upper estimate of measure and dimension of Hausdorff--Lebesgue is Lyapunov
dimension.
Recall the definition of Lyapunov dimension [2,5].

\textbf{Definition 8.} The local Lyapunov dimension of the map $F$
at the point $x$ is the number
$$
\dim_L(F,x)=j+s,
$$
where $j$ is the largest integer from interval $[1,n]$ such that
$$
\alpha_1(T_xF)\ldots\alpha_j(T_xF)\ge 1
$$
and $s$ is such that $s\in[0,1]$ and
$$
 \alpha_1(T_xF)\ldots\alpha_J(T_xF)\alpha_{j+1}(T_xF)^s= 1.
 $$
 By definition, in the case $\alpha_1(T_xF)<1$ we have $\dim_L(F,x)=0$
and in the case
 $$
 \alpha_1(T_xF)\ldots\alpha_n(T_xF)\ge 1
$$
we have $\dim_L(F,x)=n$.

\textbf{Definition 9.} The Lyapunov dimension of the map $F$ on the set $K$ is the number
$$
\dim_L(F,K)=\sup\limits_K\dim_L(F,x).
$$

\textbf{Definition 10.} A local Lyapunov dimension of the sequece
of maps $F^m$ at the point $x$ is a number
$$
\dim_Lx=\limsup\limits_{m\to+\infty}\dim_L(F^m,x).
$$

\textbf{Definition 11.} The Lyapunov dimension of maps $F^m$
on the set $K$ is a number
$$
\dim_LK=\sup\limits_K\dim_Lx.
$$

Theorem 1 implies the following result.

\textbf{Theorem 3.} Suppose that $F(K)=K$. Then $\dim_{HL}K\le\dim_LK$.

\textbf{Hypothesis.} If $F(K)=K$, then $\dim_{FHL}K\le \dim_LK$.

The theory of Lyapunov dimension of attractors is well developed [2,5--8]. For many
classical attractors the estimates and formulas of Lyapunov dimension
are obtained.
Consider such attractors.

Consider the dynamical systems generated by the differential equations
$$
\frac{dX}{dt}=f(X),\quad X\in R^n,\quad t\in R^1\eqno(8)
$$
or by the difference equations
$$
X(t+1)=f(X(t)),\quad X\in R^n,\quad t\in Z.\eqno(9)
$$
Here $Z$ is a set of integers, $f(X)$ is a vector-function: $R^n\to R^n$.
We assume that the trajectory
$X(t,X_0)$ of equation (8) is  uniquely determined for $t\in R^1$.
Here $X(0,X_0)=X_0$.

\textbf{Definition 12.} We say that $K$
is invariant if $X(t,K)=K$, $\forall t\in R^1$. Here
$$
X(t,K)=\{X(t,X_0)\,\vert\,X_0\in K\}.
$$

\textbf{Definition 13.} We say that the invariant set $K$
is locally attractive if for a certain
$\varepsilon$-neighborhood $K(\varepsilon)$ of $K$ the relation
$$
\lim\limits_{t\to+\infty}\rho(K,X(t,x_O))=0,\quad\forall x_0\in K(\varepsilon)
$$
is satisfied.

Here $\rho(K,x)$ is a distance from the point $x$ to the set $K$, defined as
$$
\rho(K,X)=\inf\limits_{Y\in K}|Y-X|,
$$
$|\cdot|$ is Euclidian norm in $R^n$,
$$
K(\varepsilon)=\{Y\,|\,\rho(K,Y)\le\varepsilon\}.
$$

\textbf{Definition 14.} We say that the invariant set $K$ is globally
attractive if
$$
\lim\limits_{t\to+\infty}\rho(K,X(t,x_0))=0,\quad\forall x_0\in R^n.
$$

\textbf{Definition 15.} We say that $K$ is

1) an attractor if it is an invariant closed and locally attractive set

2) a global attractor if it is an  invariant closed and globally attractive set.
Consider a Lorenz system [9]
$$
\dot x=-\sigma(x-y),\quad\dot y=rx-y-xz,\quad\dot z=-bz+xy,\eqno(10)
$$
where $\sigma>0$, $r>1$, $b\in[0,4]$.

\textbf{Theorem 4.} [8] If
$$
\frac{2(\sigma+b+1)}{\sigma+1+\sqrt{(\sigma-1)^2+4\sigma r}}>1,
$$
then any solution of system (8) tends to equilibrium as $t\to+\infty$. If
$$
\frac{2(\sigma+b+1)}{\sigma+1+\sqrt{(\sigma-1)^2+4\sigma r}}\le 1,
$$
then
$$
\dim_LK=3-\frac{2(\sigma+b+1)}{\sigma+1+\sqrt{(\sigma-1)^2+4\sigma r}}.\eqno(11)
$$
Here $K$ is a global attractor of system (10).

From Theorems 3 and 4 it follows that for a global attractor $K$ of system (10)
with $\sigma=10$, $b=8\!/3$, $r=28$ we have
$$
\dim_{HL}K\le 2.4014.
$$

Consider now a local attractor $K$ of system (10), which does not involve
equilibria.

It is well known that for system (10) we have
$$
\omega_3(T_xF^t)=e^{-(\sigma+b+1)t},\quad\forall x\in R^3.\eqno(12)
$$
Here $F^t$ is a shift operator along trajectories of system (10).

It is also well known that if
$$
\sup\limits_K\alpha_1(T_xF^t)\le e^{at},
$$
where $a$ is positive number, then
$$
\sup\limits_K[\alpha_1(T_xF^t)\alpha_2(T_xF^t)]\le e^{at}.
\eqno(13)
$$
It follows from the fact that in this case either the first,
either the second Lyapunov exponent is equal to zero.

Theorem 1 and [10] implies that in this case we have
$$
\dim_{HL}K\le 2+a\!/(\sigma+b+1+a).
$$
The numerical results give for
$\sigma=10$, $b=8\!/3$, $r=28$, $a\le 0.829$.

Consequently in this case ww have
$$
\dim_HK\le\dim_{HL}K\le 2.058.\eqno(14)
$$

\section{Lower estimates of Hausdorff-Lebesgue measure}

Consider one-parameter group of diffeomorphisms $F^t$, $k$-dimensional smooth
manifold $K$, $k$-dimensional segment of surface $K_1\subset K$.

\textbf{Theorem 5.} Suppose that $F^tK_1\subset K$, $\forall t>0$ $\mu_{HL}(K_1,k)>0$,
and for a certain
$t_i\to+\infty$ the following conditions
$$
\lim\limits_{t_j\to+\infty}\inf\limits_K(\alpha_k(T_xF^{t_j})
-2\alpha_{k+1}(T_xF^{t_j})>0,\eqno(15)
$$
$$
\lim\limits_{t_j\to+\infty}(\inf\limits_K\omega_k(T_xF^{t_j})=+\infty\eqno(16)
$$
are satisfied.
Then
$$
\mu_{HL}(K,k)=+\infty.
$$

Proof. It is obvious that here $\mu_{FHL}(K,k)=\mu_{HL}(K,k)$.

Relation (16) implies that for any $R>0$ there exists $\tau>0$ such
that
$$
\inf\limits_K\omega_k(T_xF^{\tau})\ge R.
$$
We choose a number  $\widetilde{\varepsilon}$ such that in any
$[\sup\limits_K\alpha_1(T_xF^\tau)]\widetilde{\varepsilon}$
--neighborhood of the point $x\in K$ there exists
a linearization proceduree
$$
F^\tau(x+h)-F^\tau(x)=(T_xF^\tau)h+o(h).
$$

Suppose that $\mu_{FHL}(K_1,k)=\mu_{HL}(K_1,k)=\nu>0$. Consider a covering $K_1$
by  disjoint cubes $C_i$ with sides $2\varepsilon<2\widetilde{\varepsilon}$ and
centers at the points $x_i$. Definition $\mu_{FHL}$ implies that the number
of these cubes is as follows
$$
N\approx\frac{\nu}{\varepsilon^k}.
$$

These cubes involve oriented cubes $\widetilde{C}_i$ with sides
$2\varepsilon\!/\sqrt n$. Obviously, $T_{x_i}F^\tau\widetilde{C}_i$ is a
parallelepiped with sides $2\varepsilon\alpha_j(T_{x_i}F^\tau)\!/\sqrt n$ and
$$
T_{x_i}F^\tau C_i\supset T_{x_i}F^\tau\widetilde{C}.
$$
This parallelepiped contains
$$
\left(\frac{\alpha_1(T_{x_i}F^\tau)}{\alpha_{k+1}(T_{x_i}F^\tau)}-1\right)\cdots
\left(\frac{\alpha_k(T_{x_i}F^\tau)}{\alpha_{k+1}(T_{x_i}F^\tau)}-1\right)
$$
disjoint cubes with sides $2\alpha_{k+1}(T_{x_i}F^\tau)\varepsilon$.
Any such cube contains points from $F^\tau K_1$.
Then from (15) it follows that
$$
\mu_{HL}(F^\tau K_1,k)\ge \left(\frac{\nu}{\varepsilon^k}\right)
\frac{1}{2^kn^{k\!/2}}\inf\limits_K
\omega_k(T_xF^\tau)\varepsilon^k\ge\frac{R\nu}{2^kn^{k\!/2}}\eqno(17)
$$

From (17) and inclusion $F^\tau K_1\subset K$ it follows that $\mu_{HL}(K,k)=+\infty$.

Theorem 5 implies that if (15) and (16) are satisfied,
then the compact $K$ cannot be bounded closed manifold.
We show numerically
that for the Lorenz system  with $\sigma=10$, $b=2\!/3$, $r=28$ the estimation $\inf\limits_K\alpha_1(T_xF^\tau)\ge{\rm{exp}}(0.788t)$ is valid.
Consequently a smooth manifold cannot be an attractor of Lorenz system.

Let us give a geometric interpretation of the proof of theorem.

\begin{figure}[thpb]
\centering
\includegraphics[scale=1.0]{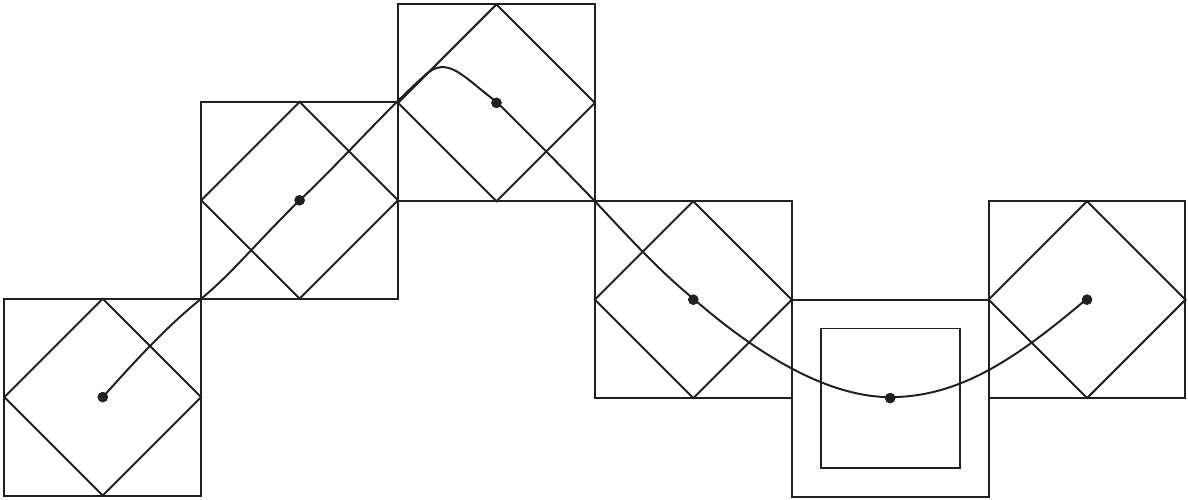}
\caption{}
\end{figure}

In Fig.1 the curve $K_1$ is covered by disjoint cubes $C_i$ with centers $x_i$
on a curve and the lengths of sides $2\varepsilon$.
Inside these cubes there are oriented cubes $\widetilde{C}_i$
with sides $\sqrt 2\varepsilon$.
The number of these cubes is $N\approx\frac{\nu}{\varepsilon}$.

\begin{figure}[thpb]
\centering
\includegraphics[scale=1.0]{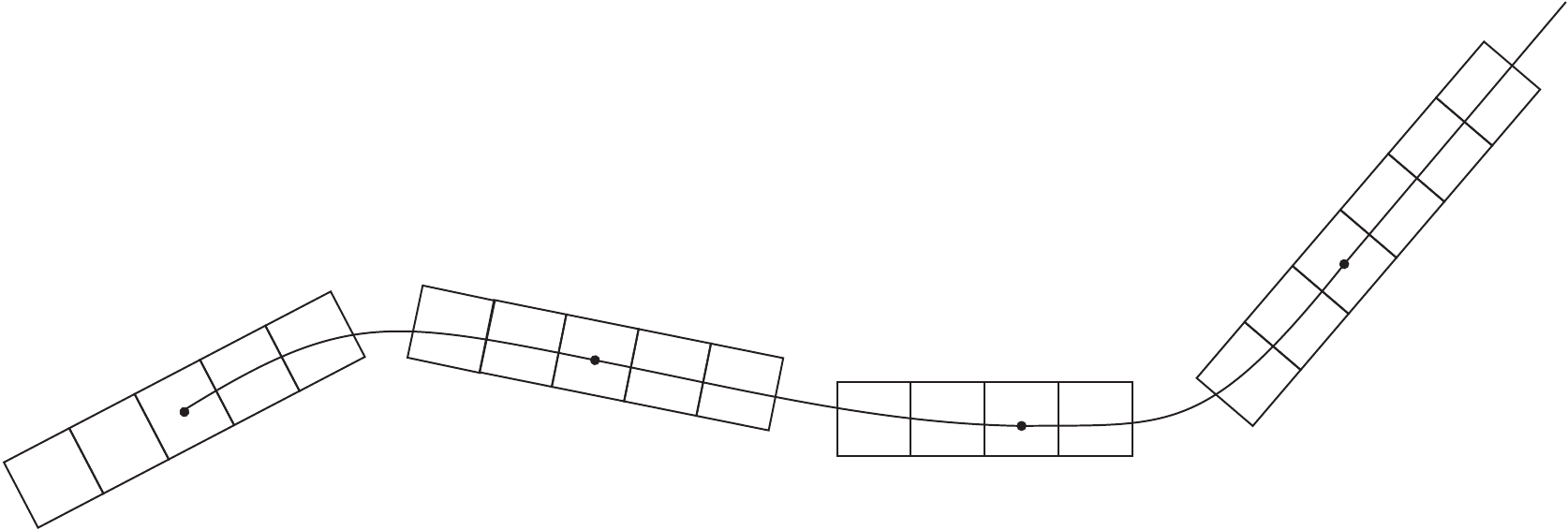}
\caption{}
\end{figure}

In Fig.2 the curve $F^\tau K_1$ is partially covered by parallelograms
$T_{x_i}F^\tau\widetilde{C}_i$.
They invove the cubes with sides $\sqrt 2\alpha_2(T_{x_i}F^\tau)\varepsilon$.
The number of all such cubes is as follows
$$
\widetilde{N}\approx\frac{\nu}{\varepsilon}\frac{\alpha_1(T_{x_i}F^\tau)}
{\sqrt 2\alpha_2(T_{x_i}F^\tau)}.
$$
It is clear that the length of curve $F^\tau K_1$ is  greater than or equal to
$$
\widetilde{N}\alpha_2(T_{x_i}F^\tau)=\nu\inf\limits_K\omega_1(T_xF^\tau)\ge\nu R\!/\sqrt 2.
$$
Similar consideration is valid for $k=2,\ldots,n-1$.

\textbf{Problem.} To extend Theorem 5 to any $d\in(1,n)$ and more
wide class of sets $K$.

\section*{References}

1. Hausdorff F. Dimension und \"{a}u{\ss}ere Ma{\ss}. Mathematische Annalen, 79, 1919, 157--179

2. Boichenko V.A.. Leonov G.A., Reitmann V. Dimension Theory for Ordinary Differential Equations.
Teubner, Wiesbaden, 2005. 440

3. Gantmacher F.R. The Theory of Matrices. Chelsea Publishing Company. New York, 1959

4. Boichenko V.A., Leonov G.A. Lyapunov Function, Lozinskii Norms and the Hausdorff Measure in the
Qualitative Theory of Differential Equations. American Mathematical Society Translations. Vol. 193,
N 2, 1999. 1--26

5. Leonov G.A. Strange Attractors and Classical Stability Theory. St. Petersburg University Press.
St. Petersburg. 2008. 160

6. Leonov G.A. On Estimations of the Hausdorff Dimension of Attractors. Vestnik St. Petersburg University,
Mathematics 24, N 3, 1991, 41--41

7. Leonov G.A. Lyapunov Dimension Formulas for Henon and Lorenz Attractors. St. Petersburg
Mathematical Journal. 13, N 3, 2002, 453--464

8. Leonov G.A. Lyapunov Dimension Formulas for Lorenz-like systems. International Journal of Bifurcation
and Chaos. 2016

9. Lorenz E.N. Deterministic Nonperiodic Flow Journal of Atmospheric Science. 20, 1963, 130--141

10. Constantin P.,  Foias C.,  Temam R.
Attractors representing turbulent flows, Mem. Am. Math. Soc. 53(314) (1985).

\end{document}